\UseRawInputEncoding
\documentclass[12pt]{article}

\usepackage{mathrsfs}
\usepackage{amssymb}
\usepackage{amsmath}
\usepackage{amsfonts}
\usepackage{amstext}
\usepackage{amsthm}
\usepackage{graphicx}
\usepackage{subfig}
\usepackage{color}
\usepackage{indentfirst,latexsym,bm}
\usepackage{mathrsfs}
\usepackage{cite}
\usepackage[all]{xy}
\usepackage{newlfont}

\usepackage{enumerate}
\newtheorem{thm}{Theorem}[section]

\newtheorem{example}[thm]{Example}
\numberwithin{equation}{section}

\newtheorem{theorem}{Theorem}[section]
\newtheorem{corollary}[subsection]{Corollary}
\newtheorem{lemma}[theorem]{Lemma}
\newtheorem{remark}[theorem]{Remark}

\newtheorem{proposition}[theorem]{Proposition}
\numberwithin{equation}{section}
\allowdisplaybreaks

\begin{document}

\title{Steady periodic hydroelastic waves and Wilton ripples with constant vorticity
\thanks{ The research of Zhang was supported by National Natural Science Foundation of China (No. 12301133), the Postdoctoral Science Foundation of China (No. 2023M741441, No. 2024T170353) and Jiangsu Education Department (No. 23KJB110007).}}
\author{Yong Zhang
\thanks{Corresponding author. School of Mathematical Sciences, Jiangsu University, Zhenjiang 212013, People's Republic of China. E-mail: zhangyong@ujs.edu.cn} \,,
Ke Zhou
\thanks{School of Mathematical Sciences, Jiangsu University, Zhenjiang 212013, People's Republic of China. E-mail: zhouke199908@163.com}
}
\date{}
\maketitle

\renewcommand{\abstractname}{Abstract}

\begin{abstract}
We study two-dimensional steady periodic hydroelastic waves in a finite-depth fluid with constant vorticity, with the free surface modeled as an extensible nonlinear hyperelastic membrane whose stored energy is $E(\nu,\mu)$. A key point is that the material coordinate of the membrane and the conformal boundary coordinate are not identified. We introduce the material reparametrization explicitly, use the tangential membrane equilibrium together with the period constraint to eliminate it locally by the implicit function theorem, and thereby obtain a fourth-order quasilinear pseudodifferential equation for the free-surface profile alone. We prove a converse reconstruction, including the material parametrization, so that the reduced fixed-domain problem is locally equivalent to the full hydroelastic free-boundary system. For fixed vorticity, isolated simple roots of the resulting dispersion relation generate local Crandall--Rabinowitz branches. Our main result concerns isolated $1{:}2$ resonances. We identify the appropriate nonlinear invariant subspace, derive the exact double-kernel conditions, and compute the quadratic resonant coefficient for the material-relaxed membrane. Under the nonvanishing of this coefficient, a Lyapunov--Schmidt secondary-bifurcation lemma adapted from Shearer yields a period-doubling secondary branch from the $2n$-mode primary family; the secondary profiles contain both the $n$- and $2n$-harmonics and have minimal period $2\pi/n$. For the quadratic constitutive law $E(\nu,\mu)=\frac A2(\nu-1)^2+\frac B2\mu^2$, $A,B>0$, the resonant coefficient is given in closed form and is strictly positive under an explicit depth condition. We finally record the exact stagnation-line criterion for the underlying laminar flows, without imposing unsupported conclusions on the topology of nonlinear critical layers.
\end{abstract}



\tableofcontents

\section{Introduction}
This paper considers two-dimensional steady periodic hydroelastic waves
propagating with constant vorticity. These waves occur on the surface of water waves, beneath a thin, frictionless elastic membrane that exhibits nonlinear responses to bending, compression, and stretching. Following Antman's formulation \cite[Chapter 4]{Antman}, we treat the elastic surface's two-dimensional cross-section as a thin, unshearable hyperelastic Cosserat rod with zero density.  This setup, termed a hydroelastic travelling wave problem, has significant physical applications, such as modeling very large floating structures or flow beneath ice sheets.

Mathematical studies of such waves began with a linear theory in the nineteenth century \cite{Greenhill}. However, detailed local and global analysis of exact nonlinear models on the surface of an infinitely deep irrotational fluid has only recently advanced. A seminal work on hydroelastic waves is due to Toland \cite{Toland1}, who established the existence for a specific class of membranes (possessing infinite elastic energy beyond fixed bending/stretching limits) by formulating the problem variationally and maximizing a Lagrangian over admissible functions. Subsequent work has generalized this theory. In \cite{Toland}, the author extended it to include membranes with positive density, subject to the convexity of the variational problem. In \cite{PlotnikovT}, Plotnikov and Toland removed this convexity restriction, addressing the general case of membranes with positive mass using Young measures to manage non-convexity. Combining calculus of variations with a Riemann-Hilbert formulation, Baldi and Toland \cite{BaldiT1} established existence of such waves without vorticity. Additionally, they performed bifurcation analysis \cite{BaldiT} to establish the existence of families of small-amplitude irrotational hydroelastic waves.  More recently, for periodic interfacial hydroelastic travelling waves in infinite depth, Akers, Ambrose and Sulon proved the existence of a global curve \cite{AkersAS} and demonstrated multiple bifurcations and ripples \cite{AkersAS1}. In fact, for specific parameter values, a two-dimensional kernel emerges across different water wave models, as demonstrated in \cite{AkersAS1,DaiXZ,EhrnstromEW}. Within this resonant regime characterized by the 2D kernel, it is natural to seek families of traveling wave solutions. These solutions feature two resonant leading-order modes and are known as Wilton ripples \cite{Wilton}. Regarding periodic hydroelastic waves in polar regions, we refer to the recent work \cite{MatiocP}. For hydroelastic solitary waves in deep water, relevant numerical studies have been conducted in \cite{GaoWM,GaoWVB,MilewshiVBW}. Furthermore, some results on time-dependent hydroelastic waves are presented in \cite{LiuA, LiuA1}.

Recent work of Martin and P\u{a}r\u{a}u \cite{MartinParau2026} established local travelling hydroelastic waves with constant vorticity and finite depth for a special Cosserat ice-sheet model, including parameter regimes with stagnation points. The contribution of the present paper is different. We retain the extensible hyperelastic membrane law of Baldi and Toland, encoded by a smooth stored energy $E(\nu,\mu)$ satisfying the rest-state assumptions, and we keep the material parametrization distinct from the conformal parametrization of the free boundary. The tangential equilibrium is first reduced to a scalar first integral and then, together with the period normalization of the material map, solved locally by the implicit function theorem. This produces a genuine one-profile formulation without suppressing the material degree of freedom. We then analyze resonant double eigenvalues in the presence of constant vorticity. At an isolated $1{:}2$ resonance we derive the quadratic interaction coefficient responsible for period-doubling and prove a secondary bifurcation theorem when this coefficient is nonzero. For the quadratic stored energy $E=\frac A2(\nu-1)^2+\frac B2\mu^2$, we obtain a closed formula and a simple sufficient depth condition implying strict nondegeneracy.

The physical model comprises a finite-depth fluid layer with a rigid bottom and an elastic surface membrane. The flow is assumed to be two-dimensional and steady, and the membrane dynamics are governed by gravity and elasticity. We adopt the physical configuration presented in \cite{BaldiT,BaldiT1}.

\subsection{The membrane elasticity}

Let $a\in\mathbb R$ denote the material coordinate on the undeformed membrane, normalized by
\[
a\mapsto a+2\pi
\]
over one reference period.  The deformed membrane is described by
\begin{equation}\label{zy}
\mathbf r_{\mathrm m}(a)=\bigl(u_{\mathrm m}(a),v_{\mathrm m}(a)\bigr),
\qquad
\mathbf r_{\mathrm m}(a+2\pi)=\mathbf r_{\mathrm m}(a)+(2\pi,0).
\end{equation}
The material stretch and bending strain are
\begin{equation}\label{eq1.1}
\nu(a):=\left|\partial_a\mathbf r_{\mathrm m}(a)\right|>0,
\qquad
\mu(a):=\partial_a\vartheta(a),
\end{equation}
where $\vartheta$ is the tangent angle.  Hence the geometric curvature is
\begin{equation}\label{eq1.2}
\sigma=\frac{\mu}{\nu}.
\end{equation}

We assume that the membrane is hyperelastic in the sense of Antman \cite[Chapter 4]{Antman}.

\medskip
\noindent
\emph{(H1) (Hyperelasticity).}
There exists a smooth stored energy
\[
E=E(\nu,\mu),\qquad \nu>0,\quad \mu\in\mathbb R,
\]
such that the elastic energy of a material segment is
\[
\mathcal E(\mathbf r_{\mathrm m})
=\int E(\nu(a),\mu(a))\,da .
\]

We use the notation
\[
E_1=\partial_\nu E,\qquad E_2=\partial_\mu E,
\qquad E_{ij}=\partial_i\partial_j E,
\]
and similarly for higher derivatives.  Unless an argument is displayed explicitly,
all constitutive derivatives used in the bifurcation coefficients are evaluated at
the rest state $(\nu,\mu)=(1,0)$.

\medskip
\noindent
\emph{(H2) (Stress-free strictly convex rest state).}
\[
E(1,0)=E_1(1,0)=E_2(1,0)=0,\qquad
E_{12}(1,0)=0,
\]
and
\[
E_{11}(1,0)>0,\qquad E_{22}(1,0)>0.
\]

For an extensible unshearable hyperelastic rod, the tangential and normal equilibrium equations are
\begin{equation}\label{eq1.tan}
\nu\,\partial_a E_1+\mu\,\partial_a E_2=0,
\end{equation}
and
\begin{equation}\label{eq1.5}
P=
\frac1\nu\partial_a\!\left(\frac{\partial_aE_2}{\nu}\right)
-\sigma E_1,
\end{equation}
respectively; see \cite{BaldiT,BaldiT1}.  Here $P$ is the normal elastic load exerted by the fluid on the membrane.  Equation \eqref{eq1.tan} is essential: it determines the material stretch (equivalently, the material reparametrization) and cannot in general be replaced by a choice of geometric parametrization.

\subsection{The free-boundary problem}

Let $\Omega$ be the fluid region bounded below by
\[
\mathcal B:=\{(X,0):X\in\mathbb R\}
\]
and above by the membrane $\mathcal S$.  In the travelling frame the stream function
$\Psi$ satisfies
\begin{equation}\label{eq1.6}
\left\{
\begin{array}{ll}
\Delta\Psi=\gamma & \text{in }\Omega,\\
\Psi=0 & \text{on }\mathcal S,\\
\Psi=-m & \text{on }\mathcal B,
\end{array}\right.
\end{equation}
where $\gamma\in\mathbb R$ is the constant vorticity and $m$ is the relative mass flux.
With the normalization of the membrane load in \eqref{eq1.5}, Bernoulli's law on the
free surface is
\begin{equation}\label{eq1.7}
\frac12|\nabla\Psi|^2+gY+P=\frac Q2
\qquad\text{on }\mathcal S,
\end{equation}
or equivalently
\begin{equation}\label{eq1.7b}
|\nabla\Psi|^2+2gY+2P=Q
\qquad\text{on }\mathcal S .
\end{equation}
Thus the full hydroelastic problem consists of the fluid equations
\eqref{eq1.6}--\eqref{eq1.7b}, the normal membrane law \eqref{eq1.5}, and the tangential
equilibrium \eqref{eq1.tan}.  The factor $2$ in front of $P$ in \eqref{eq1.7b} is simply the consequence of writing Bernoulli's law after multiplication by $2$.

\subsection{The outline of the paper}

The remainder of this paper is organized as follows. In Section 2, we introduce the conformal boundary coordinate and a separate material reparametrization.  The tangential equilibrium and the period constraint are solved locally by the implicit function theorem, after which the full free-boundary problem is reduced to a fourth-order quasilinear pseudodifferential equation for the free-surface profile; a converse reconstruction is proved as well. Section 3 establishes the Fredholm structure of the resulting linearization and constructs local branches from isolated simple roots by the Crandall--Rabinowitz theorem \cite{CrandallR}. Section 4 treats isolated $1{:}2$ resonances, computes the material-relaxed quadratic coupling coefficient, and applies a Lyapunov--Schmidt secondary-bifurcation lemma adapted from Shearer. Section 5 records the stagnation-line criterion for the laminar flows. We deliberately separate this statement from the finer nonlinear topology of critical layers, which would require an additional interior-flow expansion.

\section{Conformal reduction and elimination of the material parametrization}
\label{sec:reduction}

Let
\[
\mathcal R_h:=\{(x,y)\in\mathbb R^2:-h<y<0\}.
\]
As usual for a $2\pi$-periodic strip-like domain, we use a conformal map
$U+iV:\mathcal R_h\to\Omega$ (see following Figure \ref{fig1}) normalized by
\begin{equation}\label{eq2.1}
U(x+2\pi,y)=U(x,y)+2\pi,\qquad
V(x+2\pi,y)=V(x,y),
\end{equation}
with
\[
V(x,0)=v(x),\qquad V(x,-h)=0,\qquad [v]=h.
\]
The free surface $\mathcal{S}$ is therefore parametrized geometrically by
\begin{equation}\label{eq2.3}
\mathbf R(x)
:=\bigl(U(x,0),V(x,0)\bigr)
=\bigl(x+\mathcal C_h(v-h)(x),v(x)\bigr),
\end{equation}
up to an irrelevant horizontal translation.
\begin{figure}[ht]
\centering
\includegraphics[width=0.75\textwidth]{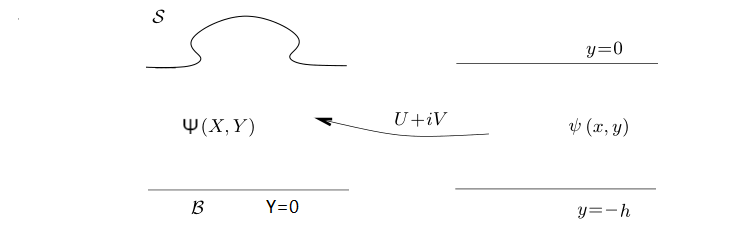}
\caption{The conformal parametrization of the fluid domain.}
\label{fig1}
\end{figure}
Here $\mathcal C_h$ is the periodic
Hilbert transform for the strip.
 In particular,
\begin{equation}\label{eq2.5}
V_y(x,0)=1+\mathcal C_h(v').
\end{equation}

Define the geometric metric factor, the tangent-angle derivative, and the curvature by
\begin{equation}\label{eq2.14}
\begin{aligned}
\Omega(v)
&:=|\mathbf R_x|
=\left((1+\mathcal C_h(v'))^2+(v')^2\right)^{1/2},\\
\mathcal N(v)
&:=v''+\mathcal C_h(v')v''-v'\mathcal C_h(v''),\\
\Theta_x(v)
&:=\frac{\mathcal N(v)}{\Omega(v)^2},
\qquad
\sigma(v):=\frac{\Theta_x(v)}{\Omega(v)}
=\frac{\mathcal N(v)}{\Omega(v)^3}.
\end{aligned}
\end{equation}

\subsection{Material coordinate versus conformal coordinate}

The variable $x$ in \eqref{eq2.3} is a conformal boundary coordinate, not a material
coordinate.  Let $a$ denote the material coordinate from Section 1 and introduce an
orientation-preserving degree-one reparametrization
\begin{equation}\label{eq2.kappa}
a=\kappa(x),\qquad
\kappa(x+2\pi)=\kappa(x)+2\pi,\qquad
\rho(x):=\kappa'(x)>0,
\end{equation}
such that
\[
\mathbf r_{\mathrm m}(\kappa(x))=\mathbf R(x).
\]
Then
\begin{equation}\label{eq2.strains}
\nu\circ\kappa=\frac{\Omega}{\rho},
\qquad
\mu\circ\kappa=\frac{\Theta_x}{\rho}
=\left(\nu\circ\kappa\right)\sigma.
\end{equation}
Thus the geometric curvature is invariant under reparametrization, while the material
strains $\nu$ and $\mu$ are not.

The tangential equilibrium \eqref{eq1.tan} has the first integral
\begin{equation}\label{eq2.firstintegral}
E(\nu,\mu)-\nu E_1(\nu,\mu)-\mu E_2(\nu,\mu)=c,
\end{equation}
where $c$ is constant along the membrane.  To solve simultaneously for the material
stretch and the material reparametrization, define
\begin{equation}\label{eq2.H}
\mathscr H(\nu,s)
:=
E(\nu,\nu s)-\nu E_1(\nu,\nu s)-\nu s E_2(\nu,\nu s).
\end{equation}
By (H2),
\[
\mathscr H(1,0)=0,\qquad
\partial_\nu\mathscr H(1,0)=-E_{11}(1,0)<0,
\qquad
\partial_s\mathscr H(1,0)=0.
\]

Let
\[
\mathcal X^{4,\alpha}_0
:=\{w\in C^{4,\alpha}_{2\pi}: [w]=0\}.
\]

\begin{lemma}[Elimination of the material reparametrization]\label{lem:materialIFT}
There exist a neighbourhood $\mathcal N$ of $(0,0)$ in $\mathbb R^2$, a neighbourhood
$\mathcal U$ of $0$ in $\mathcal X^{4,\alpha}_0$, a smooth scalar function
\[
\mathcal V:\mathcal N\longrightarrow (0,\infty),
\qquad (s,c)\longmapsto \mathcal V(s,c),
\]
and a smooth functional
\[
\mathfrak c:\mathcal U\longrightarrow\mathbb R,
\]
with
\[
\mathcal V(0,0)=1,\qquad \mathfrak c(0)=0,
\]
such that, for every $w\in\mathcal U$,
\begin{equation}\label{eq2.Vdef}
\mathscr H\bigl(\mathcal V(\sigma(w),\mathfrak c(w)),\sigma(w)\bigr)
=\mathfrak c(w)
\end{equation}
pointwise and
\begin{equation}\label{eq2.periodconstraint}
\left[
\frac{\Omega(w)}
{\mathcal V(\sigma(w),\mathfrak c(w))}
\right]=1.
\end{equation}
Consequently,
\begin{equation}\label{eq2.reducedstrains}
\widehat\nu(w)
:=\mathcal V(\sigma(w),\mathfrak c(w)),
\qquad
\widehat\mu(w):=\widehat\nu(w)\sigma(w),
\qquad
\rho(w):=\frac{\Omega(w)}{\widehat\nu(w)}
\end{equation}
determine a unique positive degree-one material density $\rho$ near the rest state.
After imposing the harmless material-phase normalization $\kappa(0)=0$, the material
reparametrization is uniquely determined by
\[
\kappa(x)=\int_0^x\rho(t)\,dt.
\]
Moreover,
\begin{equation}\label{eq2.firstvariations}
D\mathfrak c(0)=0,\qquad
D\widehat\nu(0)[f]=0,\qquad
D\widehat\mu(0)[f]=f''.
\end{equation}
\end{lemma}

\begin{proof}
First consider the finite-dimensional equation
\[
\mathscr H(\nu,s)-c=0.
\]
Since $\partial_\nu\mathscr H(1,0)=-E_{11}(1,0)\neq0$, the ordinary implicit
function theorem gives a smooth scalar function $\nu=\mathcal V(s,c)$ for $(s,c)$
near $(0,0)$.  At the rest state,
\[
\partial_c\mathcal V(0,0)=-\frac1{E_{11}},
\qquad
\partial_s\mathcal V(0,0)=0.
\]
Next define the Banach-space map
\[
\Phi:\mathcal X^{4,\alpha}_0\times\mathbb R\longrightarrow\mathbb R,
\qquad
\Phi(w,c)
:=
\left[\frac{\Omega(w)}{\mathcal V(\sigma(w),c)}\right]-1.
\]
For $w$ sufficiently small this map is smooth, because
$w\mapsto\Omega(w)$ is smooth into $C^{3,\alpha}_{2\pi}$ and
$w\mapsto\sigma(w)$ is smooth into $C^{2,\alpha}_{2\pi}$.  At $(w,c)=(0,0)$,
\[
\partial_c\Phi(0,0)
=-\partial_c\mathcal V(0,0)
=\frac1{E_{11}}>0.
\]
The Banach-space implicit function theorem therefore yields a unique smooth functional
$c=\mathfrak c(w)$ satisfying \eqref{eq2.periodconstraint}.  Since
$D\Omega(0)[f]=\mathcal C_h(f')$ has zero mean and
$D\sigma(0)[f]=f''$, differentiating $\Phi(w,\mathfrak c(w))=0$ at $w=0$ gives
$D\mathfrak c(0)=0$, and then \eqref{eq2.firstvariations} follows from the derivatives
of $\mathcal V$ above.

Finally, $\rho(w)>0$ for $w$ sufficiently small and \eqref{eq2.periodconstraint} gives
$[\rho]=1$.  Hence
\[
\kappa(x):=\int_0^x\rho(t)\,dt
\]
satisfies $\kappa(0)=0$ and \eqref{eq2.kappa}.  The density $\rho$ is unique by the two
implicit-function constructions.  Any two primitives of the same $\rho$ differ only by an
additive constant, so the normalization $\kappa(0)=0$ makes $\kappa$ unique as well.
\end{proof}

For later use write
\begin{equation}\label{eq2.Ehat}
\widehat E_j(w)
:=E_j\bigl(\widehat\nu(w),\widehat\mu(w)\bigr),
\qquad j=1,2.
\end{equation}
Transforming \eqref{eq1.5} from the material coordinate $a$ to the conformal coordinate
$x$, using $\partial_a=\rho^{-1}\partial_x$ and
$\nu=\Omega/\rho$, gives the exact identity
\begin{equation}\label{eq2.15}
P(w)
=
\frac1{\Omega}
\left(\frac{\widehat E_2'(w)}{\Omega}\right)'
-\sigma(w)\widehat E_1(w).
\end{equation}
Indeed, the factors $\rho$ cancel:
\[
\frac1\nu\partial_a
\left(\frac{\partial_aE_2}{\nu}\right)
=
\frac1\Omega\partial_x
\left(\frac{\partial_xE_2}{\Omega}\right).
\]
This cancellation is the reason why the material variable can be eliminated without
increasing the differential order of the free-surface equation.

\subsection{Fluid part and the scalar profile equation}

Let
\[
\psi(x,y):=\Psi(U(x,y),V(x,y)).
\]
Then
\[
\eta:=\psi-\frac{\gamma}{2}V^2+m
\]
is harmonic in $\mathcal R_h$ with
\[
\eta(x,-h)=0,\qquad
\eta(x,0)=m-\frac{\gamma}{2}v(x)^2.
\]
The strip Dirichlet--Neumann formula gives
\begin{equation}\label{eq2.13}
V_y(x,0)=1+\mathcal C_h(v'),
\qquad
\eta_y(x,0)
=
\frac mh-\frac{\gamma[v^2]}{2h}
-\gamma\mathcal C_h(vv').
\end{equation}
By conformality,
\[
|\nabla\Psi|^2
=
\frac{(\eta_y+\gamma vV_y)^2}{\Omega^2}
\qquad\text{on }y=0.
\]
Hence \eqref{eq1.7b} becomes
\begin{equation}\label{eq2.12}
\frac{(\eta_y+\gamma vV_y)^2}{\Omega^2}
+2gv+2P(w)=Q.
\end{equation}

Substituting \eqref{eq2.13} and \eqref{eq2.15} into \eqref{eq2.12}, and multiplying
by $\Omega^3$, yields the reduced scalar equation
\begin{eqnarray}\label{eq2.16}
&&
\left(
\frac mh-\frac{\gamma[v^2]}{2h}
-\gamma\mathcal C_h(vv')
+\gamma v(1+\mathcal C_h(v'))
\right)^2\Omega
\nonumber\\
&&\qquad
=(Q-2gv)\Omega^3
-2\widehat E_2''\,\Omega
+2\widehat E_2'\,\Omega'
+2\mathcal N(v)\widehat E_1 .
\end{eqnarray}

The admissible profiles additionally satisfy
\[
[v]=h,\qquad v>0,
\]
the boundary map in \eqref{eq2.3} is injective, and $\Omega$ is bounded away from zero.

\begin{theorem}[Local equivalence with the full hydroelastic problem]\label{thm2equiv}
Fix $h>0$ and $\alpha\in(0,1)$.  There exists a neighbourhood $\mathcal U$ of
the flat profile $v\equiv h$ in $C^{4,\alpha}_{2\pi}$ with the following property.

If a solution of the full physical problem
\eqref{eq1.tan}, \eqref{eq1.5}, \eqref{eq1.6}, \eqref{eq1.7b}
has free surface in $\mathcal U$ and a degree-one material parametrization close to the
identity, then its conformal profile satisfies \eqref{eq2.16} with the reduced strains
\eqref{eq2.reducedstrains}.

Conversely, if $v\in\mathcal U$ satisfies \eqref{eq2.16}, $[v]=h$, and the geometric
injectivity and nondegeneracy conditions above, then \eqref{eq2.reducedstrains}
constructs a positive $\rho$ and hence a degree-one material map $\kappa$.  Together
with the conformal reconstruction of the fluid, these quantities determine a classical
solution of the full system
\eqref{eq1.tan}, \eqref{eq1.5}, \eqref{eq1.6}, \eqref{eq1.7b}.
\end{theorem}

\begin{proof}
For the forward implication, the tangential balance gives
\eqref{eq2.firstintegral}.  Relations \eqref{eq2.strains} and $[\rho]=1$ then imply,
by the uniqueness in Lemma \ref{lem:materialIFT}, that the physical material strains
coincide with \eqref{eq2.reducedstrains}.  Formula \eqref{eq2.15} follows from the
change of material coordinate, and the calculation above gives \eqref{eq2.16}.

Conversely, Lemma \ref{lem:materialIFT} gives $\rho$ and $\kappa$, and therefore a
material parametrization
\[
\mathbf r_{\mathrm m}(a)
:=\mathbf R(\kappa^{-1}(a)).
\]
By construction \eqref{eq2.firstintegral} holds, hence its derivative is precisely
the tangential equilibrium \eqref{eq1.tan}; formula \eqref{eq2.15} is the normal
balance \eqref{eq1.5} written in conformal coordinates.

To reconstruct the fluid, let $V$ be the harmonic extension of $v$ with $V=0$ on
$y=-h$, let $U$ be its harmonic conjugate, and define
\[
\psi=\eta+\frac{\gamma}{2}V^2-m.
\]
The standard strip reconstruction argument
\cite[Appendix A]{ConstantinV} gives a conformal homeomorphism
$U+iV:\overline{\mathcal R_h}\to\overline\Omega$ for profiles in $\mathcal U$ satisfying
the geometric conditions.  Since
\[
\Delta\psi=\gamma|\nabla V|^2,
\]
the function $\Psi$ defined by $\Psi(U,V)=\psi$ satisfies $\Delta\Psi=\gamma$ and the
Dirichlet conditions in \eqref{eq1.6}.  Finally \eqref{eq2.16} is exactly the pullback
of \eqref{eq1.7b}.  This proves the equivalence.
\end{proof}

\begin{remark}
The reduction is local near the stress-free state.  This is exactly the regime needed
for the bifurcation analysis below.  The material coordinate has not been identified
with the conformal coordinate; rather, it has been solved for and eliminated.
\end{remark}

\section{Bifurcation from a simple eigenvalue}

Write
\[
w=v-h,\qquad [w]=0,
\]
and set
\begin{equation}\label{eq3.3}
\lambda:=\frac mh+\frac{\gamma h}{2},
\qquad
\theta:=Q-2gh-\lambda^2.
\end{equation}
At $w=0$, Lemma \ref{lem:materialIFT} gives
\[
\widehat\nu=1,\qquad \widehat\mu=0,\qquad \rho=1,\qquad P=0,
\]
so the flat laminar family is characterized by $\theta=0$.  Its stream function is
\begin{equation}\label{eq3.5}
\psi(X,Y)
=
\frac{\gamma}{2}Y^2
+\left(\frac mh-\frac{\gamma h}{2}\right)Y-m.
\end{equation}

Define
\begin{equation}\label{eq3.R}
\mathscr R(w)
:=
\frac{[w^2]}{2h}-w
+\mathcal C_h(ww')-w\mathcal C_h(w').
\end{equation}
Then \eqref{eq2.16} is equivalent to
\begin{eqnarray}\label{eq3.4}
F(\lambda,\gamma,(\theta,w))
&:=&
\lambda^2\Omega
+\gamma^2\mathscr R(w)^2\Omega
-2\lambda\gamma\mathscr R(w)\Omega
-\bigl(\theta+\lambda^2-2gw\bigr)\Omega^3
\nonumber\\
&&
+2\widehat E_2''\,\Omega
-2\widehat E_2'\,\Omega'
-2\mathcal N(w)\widehat E_1
=0.
\end{eqnarray}
Here $\Omega,\mathcal N$ are the geometric quantities in \eqref{eq2.14}, with
$v=h+w$, and $\widehat E_i$ are defined by \eqref{eq2.Ehat} through the eliminated
material parametrization.

For $\alpha\in(0,1)$ put
\[
C^{4,\alpha}_{2\pi,e,0}
:=
\{f\in C^{4,\alpha}_{2\pi}:f(-x)=f(x),\ [f]=0\},
\]
\[
X:=\mathbb R\times C^{4,\alpha}_{2\pi,e,0},
\qquad
Y:=C^{0,\alpha}_{2\pi,e}.
\]
On a sufficiently small open neighbourhood $\mathcal O\subset X$ of $(0,0)$ for which
the conformal and material metric factors are positive, the map
\[
F:\mathbb R^2\times\mathcal O\to Y
\]
is smooth and
\[
F(\lambda,\gamma,(0,0))=0.
\]

\subsection{Linearized spectrum}

Let
\[
B:=E_{22}(1,0)>0.
\]
By Lemma \ref{lem:materialIFT},
\begin{equation}\label{eq3.6}
D\Omega(0)[f]=\mathcal C_h(f'),
\qquad
D\widehat\nu(0)[f]=0,
\qquad
D\widehat\mu(0)[f]=f''.
\end{equation}
Using $E_{12}(1,0)=E_1(1,0)=0$, the flat-state linearization is therefore
\begin{equation}\label{eq3.7}
\begin{aligned}
L_{\lambda,\gamma}(\zeta,f)
&:=F_{(\theta,w)}(\lambda,\gamma,(0,0))(\zeta,f)\\
&=-\zeta
-2\lambda^2\mathcal C_h(f')
+2\lambda\gamma f+2gf
+2Bf''''.
\end{aligned}
\end{equation}
For
\[
f=\sum_{k\ge1}a_k\cos(kx)
\]
this becomes
\begin{equation}\label{eq3.9}
L_{\lambda,\gamma}(\zeta,f)
=
-\zeta-\sum_{k\ge1}D_k(\lambda,\gamma)a_k\cos(kx),
\end{equation}
where
\begin{equation}\label{eq3.10}
D_k(\lambda,\gamma)
=
2k\coth(kh)\lambda^2
-2\lambda\gamma
-\bigl(2g+2Bk^4\bigr).
\end{equation}

\begin{lemma}[Fredholm structure]\label{lem3fredholm}
For every $(\lambda,\gamma)\in\mathbb R^2$,
\[
L_{\lambda,\gamma}:
\mathbb R\times C^{4,\alpha}_{2\pi,e,0}
\longrightarrow
C^{0,\alpha}_{2\pi,e}.
\]
is Fredholm of index zero.  Its kernel is generated by the Fourier modes associated with zeros of the dispersion multiplier; more precisely, by $\cos(kx)$ whenever
$D_k(\lambda,\gamma)=0$.
\end{lemma}

\begin{proof}
Set
\[
L_0(\zeta,f):=-\zeta+2B\partial_x^4f.
\]
The operator $L_0$ is an isomorphism from
$\mathbb R\times C^{4,\alpha}_{2\pi,e,0}$ onto $C^{0,\alpha}_{2\pi,e}$.
The difference $L_{\lambda,\gamma}-L_0$ is of order at most one and is compact into
$C^{0,\alpha}$ by the compact H\"older embedding.  Thus the index is zero, and the
Fourier representation \eqref{eq3.9} gives the kernel and range mode by mode.
\end{proof}

Put
\[
T_k:=\frac{\tanh(kh)}{k}.
\]
For fixed $\gamma$, the roots of $D_k(\lambda,\gamma)=0$ are
\begin{equation}\label{eq3.11}
\lambda_{k,\pm}^*
=
\frac{\gamma T_k}{2}
\pm
\sqrt{
\frac{\gamma^2T_k^2}{4}
+\bigl(g+Bk^4\bigr)T_k
}.
\end{equation}
In particular
\[
\lambda_{k,+}^*>0,\qquad \lambda_{k,-}^*<0.
\]

Fix $n\ge1$, one sign, and put $\lambda_*=\lambda_{n,\pm}^*$.  Assume the simple-root
condition
\begin{equation}\label{eq3.simple}
D_k(\lambda_*,\gamma)\neq0
\qquad(k\neq n).
\end{equation}
Then
\[
\ker L_{\lambda_*,\gamma}
=\operatorname{span}\{(0,\cos(nx))\},
\]
and
\[
Y
=
\mathcal R(L_{\lambda_*,\gamma})
\oplus\operatorname{span}\{\cos(nx)\}.
\]
Moreover
\begin{equation}\label{eq3.trans}
F_{\lambda(\theta,w)}(\lambda_*,\gamma,0)[1,(0,\cos nx)]
=
2\left(\gamma-\frac{2\lambda_*}{T_n}\right)\cos(nx),
\end{equation}
which is not in the range because
\[
\frac{2\lambda_*}{T_n}-\gamma
=
\pm\frac2{T_n}
\sqrt{
\frac{\gamma^2T_n^2}{4}
+\bigl(g+Bn^4\bigr)T_n
}
\neq0.
\]

\begin{theorem}[Local bifurcation from a simple mode]\label{thm3.1}
Assume \eqref{eq3.simple}.  There exist $\varepsilon>0$ and a $C^1$ curve
\[
\mathcal K_{n,\pm}
=
\{(\lambda_{n,\pm}(s),(\theta(s),w(s))):|s|<\varepsilon\}
\]
of solutions of \eqref{eq3.4} such that
\[
(\lambda_{n,\pm}(0),\theta(0),w(0))
=(\lambda_{n,\pm}^*,0,0),
\]
and
\[
w(s)=s\cos(nx)+o(s)
\qquad\text{in }C^{4,\alpha}.
\]
All nontrivial solutions sufficiently close to the bifurcation point lie on this curve;
the remaining nearby solutions are those on the trivial laminar branch.
\end{theorem}

\begin{proof}
Since Lemma \ref{lem3fredholm}, \eqref{eq3.simple}, and \eqref{eq3.trans} verify the
Crandall--Rabinowitz hypotheses in Theorem \ref{thm6.1}, we apply the Crandall--Rabinowitz local bifurcation Theorem \ref{thm6.1} to finish the proof.
\end{proof}

\begin{corollary}[Physical small-amplitude waves]\label{cor3physical}
For $|s|$ sufficiently small, every nontrivial solution on $\mathcal K_{n,\pm}$
reconstructs, through Theorem \ref{thm2equiv}, a classical hydroelastic wave together
with a positive degree-one material parametrization.
\end{corollary}

\begin{proof}
Since $w(s)\to0$ in $C^{4,\alpha}$,
\[
\|\mathcal C_h(w'(s))\|_{C^0}
\le C\|w(s)\|_{C^{4,\alpha}}\to0.
\]
Thus $v=h+w>0$, the horizontal derivative
$1+\mathcal C_h(w')$ is positive, and the conformal metric factor is bounded away from
zero.  Lemma \ref{lem:materialIFT} also gives $\rho(w)>0$ for $w$ small.
Theorem \ref{thm2equiv} applies.
\end{proof}

\section{Secondary bifurcation at an isolated $1{:}2$ resonance}
\label{sec:secondary}

We now vary both $\lambda$ and $\gamma$.  Fix $n\ge1$ and work in the closed
$n$-fold symmetric spaces
\[
H^s_{n,e,0}
:=
\left\{
w=\sum_{j\ge1}a_j\cos(jnx):
\sum_{j\ge1}(1+j^2n^2)^s a_j^2<\infty
\right\},
\]
\[
H^s_{n,e}
:=
\mathbb R\oplus H^s_{n,e,0},
\qquad
X:=\mathbb R\times H^5_{n,e,0},
\qquad
Y:=H^1_{n,e}.
\]
Let $\mathcal O\subset X$ be a sufficiently small neighbourhood of the origin on which
the conformal metric and the material density $\rho(w)$ are positive.  The reduced
operator \eqref{eq3.4} defines a smooth map
\[
F:\mathbb R^2\times\mathcal O\to Y.
\]

\begin{lemma}[Smoothness and Fredholm structure]\label{lem4smooth}
The reduced map $F$ is $C^\infty$ from $\mathbb R^2\times\mathcal O$ to $H^1_{n,e}$.
For every $(\lambda,\gamma)$ its flat-state linearization
\[
L_{\lambda,\gamma}:
\mathbb R\times H^5_{n,e,0}\to H^1_{n,e}
\]
is Fredholm of index zero.
\end{lemma}

\begin{proof}
For $w\in H^5$, one has $\Omega-1\in H^4$ and $\sigma\in H^3$.
The implicit-function construction in Lemma \ref{lem:materialIFT} is smooth in Sobolev
spaces as well: the scalar constraint is smooth, $E_{11}>0$ gives the required
nondegeneracy, and hence
\[
w\mapsto
\widehat\nu(w),\widehat\mu(w),\rho(w),\widehat E_j(w)
\]
is smooth with the regularity needed in \eqref{eq3.4}.  The highest derivative is
$2B\partial_x^4$, so the same compact-perturbation argument as in Lemma
\ref{lem3fredholm}, now with
\[
L_0(\zeta,f)=-\zeta+2B\partial_x^4f,
\]
gives Fredholm index zero.
\end{proof}

The dispersion multiplier remains
\begin{equation}\label{eq4.disp}
D_k(\lambda,\gamma)
=
\frac{2\lambda^2}{T_k}
-2\lambda\gamma
-\bigl(2g+2Bk^4\bigr).
\end{equation}

\subsection{The resonant pair}

\begin{lemma}[Explicit $1{:}2$ resonance]\label{lem4respair}
The simultaneous equations
\[
D_n(\lambda_*,\gamma_*)
=
D_{2n}(\lambda_*,\gamma_*)=0
\]
are equivalent to
\begin{equation}\label{eq4.respair}
\lambda_*^2
=
15Bn^3c,
\qquad
2\lambda_*\gamma_*
=
30Bn^4c^2-2g-2Bn^4,
\qquad
c:=\coth(nh).
\end{equation}
Thus the two resonant parameter pairs are related by
$(\lambda_*,\gamma_*)\mapsto(-\lambda_*,-\gamma_*)$.
\end{lemma}

\begin{proof}
Subtracting the $n$-equation from the $2n$-equation gives
\[
2\lambda_*^2n\bigl(2\coth(2nh)-\coth(nh)\bigr)
=
30Bn^4.
\]
Since
\[
2\coth(2z)-\coth z=\tanh z=\frac1{\coth z},
\]
the first formula follows.  The second is then obtained from $D_n=0$.
\end{proof}

We impose the isolated-resonance condition
\begin{equation}\label{eq4.isolated}
D_k(\lambda_*,\gamma_*)\neq0
\qquad
(k\notin\{n,2n\}).
\end{equation}
This stronger full-period condition in particular excludes every additional kernel mode
in the $n$-fold symmetric space.

\begin{lemma}[Isolated double kernel]\label{lem4double}
Under \eqref{eq4.respair}--\eqref{eq4.isolated},
\[
\ker L_{\lambda_*,\gamma_*}
=
\operatorname{span}\{x_1,x_2\},
\qquad
x_1:=(0,\cos(2nx)),\quad
x_2:=(0,\cos(nx)).
\]
The range has codimension two.
\end{lemma}

\subsection{Invariant splitting and parameter transversality}

Define
\[
X_1
:=
\left\{
(\theta,w)\in X:
w=\sum_{j\ge1}a_j\cos(2jnx)
\right\},
\]
and let $X_2$ be its Hilbert orthogonal complement in $X$, with zero
$\theta$-component.  Define $Y_1$ analogously using the constant mode and the
$2jn$-modes, and let $Y_2=Y_1^\perp$.  Then
\[
X=X_1\oplus X_2,\qquad Y=Y_1\oplus Y_2,
\]
and, by translation invariance of all geometric and material-reduction operations,
\begin{equation}\label{eq4.invariance}
F(\lambda,\gamma,X_1\cap\mathcal O)\subset Y_1.
\end{equation}

Define
\[
\mathfrak L_1(f)
=
\frac1\pi\int_{-\pi}^{\pi}f(x)\cos(2nx)\,dx,
\qquad
\mathfrak L_2(f)
=
\frac1\pi\int_{-\pi}^{\pi}f(x)\cos(nx)\,dx.
\]
Then
\begin{equation}\label{eq4.paramder}
\begin{aligned}
F_{\lambda x}[1,x_1]
&=\left(-\frac{4\lambda_*}{T_{2n}}+2\gamma_*\right)\cos(2nx),\\
F_{\lambda x}[1,x_2]
&=\left(-\frac{4\lambda_*}{T_n}+2\gamma_*\right)\cos(nx),\\
F_{\gamma x}[1,x_1]&=2\lambda_*\cos(2nx),\\
F_{\gamma x}[1,x_2]&=2\lambda_*\cos(nx).
\end{aligned}
\end{equation}
The two $\lambda$-transversality coefficients are nonzero, and
\begin{equation}\label{eq4.detparam}
\det
\begin{pmatrix}
\mathfrak L_1(F_{\lambda x}[1,x_1])&
\mathfrak L_1(F_{\gamma x}[1,x_1])\\
\mathfrak L_2(F_{\lambda x}[1,x_2])&
\mathfrak L_2(F_{\gamma x}[1,x_2])
\end{pmatrix}
=
8\lambda_*^2
\left(\frac1{T_n}-\frac1{T_{2n}}\right)
\neq0.
\end{equation}

\subsection{Quadratic coefficient after material relaxation}

Write
\[
\ell(w):=\mathcal C_h(w').
\]
At the flat state,
\begin{equation}\label{eq4.geomexp}
\Omega
=
1+\ell(w)+\frac12(w')^2+O(\|w\|^3),
\end{equation}
and
\begin{equation}\label{eq4.sigmaexp}
\sigma
=
w''-2\ell(w)w''-w'\mathcal C_h(w'')
+O(\|w\|^3).
\end{equation}
The material reduction has no linear stretch variation:
\[
D\widehat\nu(0)=0,\qquad
D\widehat\mu(0)[f]=f''.
\]
Moreover, using (H2),
\begin{equation}\label{eq4.E2exp}
D\widehat E_2(0)[f]=Bf'',
\end{equation}
and for two directions $a,b$,
\begin{equation}\label{eq4.E2second}
D^2\widehat E_2(0)[a,b]
=
B\,D^2\sigma(0)[a,b]
+E_{222}\,a''b''.
\end{equation}
Terms involving $D^2\widehat\nu$ disappear from \eqref{eq4.E2second} because
$E_{12}(1,0)=0$.  Also $\widehat E_1=O(\|w\|^2)$, so
$\mathcal N(w)\widehat E_1(w)$ has no quadratic contribution.

From \eqref{eq4.sigmaexp},
\begin{equation}\label{eq4.sigmasecond}
\begin{aligned}
D^2\sigma(0)[a,b]
={}&-2\ell(a)b''-2\ell(b)a''\\
&-a'\mathcal C_h(b'')-b'\mathcal C_h(a'').
\end{aligned}
\end{equation}

\begin{lemma}[Material-relaxed $1{:}2$ coupling coefficient]\label{lem4q}
Let
\[
c_1:=\coth(nh),\qquad c_2:=\coth(2nh),
\]
and define
\begin{equation}\label{eq4.qdef}
\mathfrak q_n
:=
\mathfrak L_2
\left(
F_{xx}(\lambda_*,\gamma_*,0)[x_1,x_2]
\right).
\end{equation}
Then
\begin{align}\label{eq4.qexplicit}
\mathfrak q_n={}&
\gamma_*^2
+n\lambda_*\gamma_*(c_1+4c_2)
+3gnc_1+6gnc_2
-n^2\lambda_*^2(2+6c_1c_2)
\nonumber\\
&\quad
+18Bn^5c_1+6Bn^5c_2
-4E_{222}n^6 .
\end{align}
\end{lemma}

\begin{proof}
Let
\[
a=\cos(2nx),\qquad b=\cos(nx).
\]
The fluid part is unchanged by the material reduction and gives
\[
\gamma_*^2
+n\lambda_*\gamma_*(c_1+4c_2)
+3gnc_1+6gnc_2
-n^2\lambda_*^2(2+6c_1c_2).
\]

For the elastic part, \eqref{eq3.4} gives
\[
F_{\mathrm{el}}(w)
=
2\left(
\Omega\,\widehat E_2''
-\Omega'\widehat E_2'
-\mathcal N(w)\widehat E_1
\right).
\]
Set $S_{ab}:=D^2\sigma(0)[a,b]$.  By
\eqref{eq4.E2second}--\eqref{eq4.sigmasecond},
\begin{align*}
F_{\mathrm{el},ww}(0)[a,b]
={}&
2B S_{ab}''
+2E_{222}(a''b'')''
+2B\bigl(\ell(a)b^{(4)}+\ell(b)a^{(4)}\bigr)\\
&-2B\bigl(\ell(a)'b'''+\ell(b)'a'''\bigr).
\end{align*}
Using
\[
\ell(a)=2nc_2\cos(2nx),\qquad
\ell(b)=nc_1\cos(nx),
\]
and the product-to-sum identities, its $\cos(nx)$ coefficient is
\[
18Bn^5c_1+6Bn^5c_2-4E_{222}n^6.
\]
Combining the fluid and elastic contributions proves \eqref{eq4.qexplicit}.
\end{proof}

We assume the explicit nondegeneracy condition
\begin{equation}\label{eq414}
\mathfrak q_n\neq0.
\end{equation}

\begin{theorem}[Secondary $1{:}2$ bifurcation]\label{thm4.1}
Assume (H1)--(H2), \eqref{eq4.isolated} and \eqref{eq414} hold.  Let
$(\lambda_*,\gamma_*)$ be one of the resonant pairs in \eqref{eq4.respair}.
Then there exists $\varepsilon>0$ such that:

\begin{enumerate}[(i)]
\item for every $\gamma$ with $|\gamma-\gamma_*|<\varepsilon$, there is a smooth
primary branch $\mathcal P_\gamma\subset X_1$ bifurcating in the direction
$x_1=(0,\cos(2nx))$;

\item for every $\gamma=\gamma_*+\delta$ with $0<|\delta|<\varepsilon$, a smooth
secondary branch $\mathcal S_\gamma$ bifurcates from a nontrivial point of
$\mathcal P_\gamma$ in the direction $x_2=(0,\cos(nx))$;

\item every secondary solution with nonzero secondary amplitude has minimal period
$2\pi/n$.
\end{enumerate}
\end{theorem}

\begin{proof}
We apply Theorem \ref{thm6.2} with $\beta=\gamma$.  Lemmas
\ref{lem4smooth} and \ref{lem4double} give the Fredholm splitting and the
two-dimensional kernel.  The invariant-space hypothesis is
\eqref{eq4.invariance}.  The parameter determinant is \eqref{eq4.detparam}.

The interaction of two $2n$-modes contains only the $0$- and $4n$-frequencies, hence
\[
\mathfrak L_1(F_{xx}[x_1,x_1])=0.
\]
Lemma \ref{lem4q} gives
\[
\mathfrak L_2(F_{xx}[x_1,x_2])=\mathfrak q_n.
\]
Since $\mathfrak L_1(F_{xx}[x_1,x_1])=0$, the quadratic Jacobian
\eqref{eq6.quaddet} reduces in the present problem to
\[
\mathfrak L_1(F_{\lambda x}[1,x_1])\,\mathfrak q_n\neq0.
\]
Theorem \ref{thm6.2} therefore applies.

Finally, the whole reduction is carried out in the $n$-fold space
$H^5_{n,e,0}$, so every secondary solution has period $2\pi/n$.  Its nonzero
$n$-Fourier coefficient excludes every strictly smaller period.  Thus the minimal period is
exactly $2\pi/n$.
\end{proof}

\subsection{Example: Quadratic elasticity}

\begin{example}[Explicit nondegeneracy for quadratic elasticity]\label{cor4quadratic}
Assume
\begin{equation}\label{eq4.quadE}
E(\nu,\mu)
=
\frac A2(\nu-1)^2+\frac B2\mu^2,
\qquad A>0,\quad B>0.
\end{equation}
Let $(\lambda_*,\gamma_*)$ satisfy \eqref{eq4.respair} and put $c=\coth(nh)$.
Then
\begin{equation}\label{eq4.qquad}
\mathfrak q_n
=
\frac{\mathcal N_n(c)}{15Bcn^3},
\end{equation}
where
\begin{align}\label{eq4.Nquad}
\mathcal N_n(c)
={}&
B^2n^8\bigl(225c^4-435c^2+16\bigr)
\nonumber\\
&+Bgn^4(15c^2+17)+g^2.
\end{align}
In particular, if
\begin{equation}\label{eq4.depthcond}
\coth^2(nh)
\ge
\frac{29+\sqrt{777}}{30},
\end{equation}
then $\mathfrak q_n>0$, so \eqref{eq414} holds automatically.
\end{example}

\begin{proof}
For \eqref{eq4.quadE}, $E_{222}=0$.  Substitute
\eqref{eq4.respair} and
\[
c_2=\frac12\left(c+\frac1c\right)
\]
into \eqref{eq4.qexplicit}.  This gives
\eqref{eq4.qquad}--\eqref{eq4.Nquad}.  The quadratic polynomial
\[
225z^2-435z+16
\]
is nonnegative for
\[
z\ge\frac{29+\sqrt{777}}{30}.
\]
The remaining terms in \eqref{eq4.Nquad} are strictly positive.
\end{proof}

\begin{remark}
The constant $A=E_{11}$ is essential for the elimination of the material parametrization:
it is exactly the nondegeneracy $\partial_\nu\mathscr H(1,0)=-A$.  It does not appear
in the projected quadratic coefficient \eqref{eq4.qquad}; this cancellation is a
consequence of $D\widehat\nu(0)=0$ and $E_{12}(1,0)=0$.
\end{remark}

\begin{remark}[Openness of the nondegeneracy condition]
The construction of Lemma \ref{lem:materialIFT} and the coefficient
\eqref{eq4.qexplicit} depend continuously on the constitutive law in a sufficiently
small $C^3$ neighbourhood of the rest state.  Hence every quadratic law satisfying
the strict inequality in Corollary \ref{cor4quadratic} is contained in an open class
of genuinely nonlinear smooth stored energies for which the secondary branch persists.
\end{remark}

\section{Laminar stagnation lines}
The bifurcation results above concern the free surface and do not require an assertion about the topology of interior streamlines. We nevertheless record the exact stagnation criterion for the laminar flows from which the branches emerge.

For the flat state, \eqref{eq3.5} and \eqref{eq3.3} give
\begin{equation}\label{eq5.1}
(\psi_Y,-\psi_X)=\bigl((Y-h)\gamma+\lambda,0\bigr),
\qquad 0\le Y\le h.
\end{equation}
Hence, when $\gamma\neq0$, an interior horizontal line of stagnation points exists if and only if
\begin{equation}\label{eq5.2}
0<\frac{\lambda}{\gamma}<h.
\end{equation}
In that case the stagnation line is
\begin{equation}\label{eq5.3}
Y_0=h-\frac{\lambda}{\gamma}\in(0,h).
\end{equation}
For $\gamma=0$, the laminar velocity is constant and no interior stagnation line exists unless $\lambda=0$, which is not a bifurcation value in the present setting.

\begin{proposition}\label{prop5laminar}
Let $(\lambda_*,\gamma_*)$ be either a simple bifurcation point from Section 3 or an isolated resonant point from Section 4. The corresponding laminar flow possesses an interior stagnation line if and only if \eqref{eq5.2} holds with $(\lambda,\gamma)=(\lambda_*,\gamma_*)$.
\end{proposition}

\begin{proof}
This is immediate from the affine formula \eqref{eq5.1} for the horizontal relative velocity.
\end{proof}

\begin{remark}
A laminar stagnation line is highly degenerate: every point on the line is critical. Showing that a nonlaminar perturbation breaks this line into centers and saddles and creates a Kelvin cat's-eye requires a separate expansion of the interior stream function and a nondegeneracy analysis of its critical points. Such a conclusion is available in several gravity-wave settings and, for a special hydroelastic Cosserat model, in \cite{MartinParau2026}; however, it does not follow solely from the free-surface bifurcation theorem proved here. We therefore do not claim that every solution on the primary or secondary branch has exactly one critical layer.
\end{remark}

\section{Appendix: bifurcation theorems used in the paper}
For convenience we record the two local bifurcation results used above. The first is the classical theorem of Crandall and Rabinowitz \cite{CrandallR}. The second is a Hilbert-space Lyapunov--Schmidt lemma adapted from Shearer's secondary-bifurcation argument \cite{Shearer} and tailored to the invariant splitting used in Section \ref{sec:secondary}.

\begin{theorem}[Crandall--Rabinowitz]\label{thm6.1}
Let $X,Y$ be Banach spaces and let $F:I\times X\to Y$ be $C^k$, $k\ge2$, where $I\subset\mathbb R$ is open. Assume
\begin{enumerate}[(i)]
\item $F(\lambda,0)=0$ for every $\lambda\in I$;
\item for some $\lambda_*\in I$, the operator $L:=F_x(\lambda_*,0)$ is Fredholm of index zero with
\[
\ker L=\operatorname{span}\{x_*\},\qquad \operatorname{codim}\mathcal R(L)=1;
\]
\item the transversality condition
\[
F_{\lambda x}(\lambda_*,0)[1,x_*]\notin\mathcal R(L)
\]
holds.
\end{enumerate}
Then there exist $\varepsilon>0$ and $C^{k-1}$ functions $\Lambda:(-\varepsilon,\varepsilon)\to I$ and $\psi:(-\varepsilon,\varepsilon)\to X$ such that
\[
\Lambda(0)=\lambda_*,\qquad \psi(0)=x_*,
\]
and all nontrivial solutions sufficiently close to $(\lambda_*,0)$ lie on the curve
\[
(\lambda,x)=\bigl(\Lambda(s),s\psi(s)\bigr),\qquad |s|<\varepsilon.
\]
\end{theorem}

\begin{theorem}[Two-parameter Lyapunov--Schmidt lemma (after Shearer)]\label{thm6.2}
Let $X,Y$ be real Hilbert spaces and let
$F:\mathbb R^2\times X\to Y$ be $C^\infty$.  Suppose there are closed topological
decompositions
\[
X=X_1\oplus X_2,\qquad Y=Y_1\oplus Y_2,
\]
and a point $(\lambda_*,\beta_*)$ such that:

\begin{enumerate}[(B1)]
\item $F(\lambda,\beta,0)=0$ near $(\lambda_*,\beta_*)$ and
\[
F(\lambda,\beta,X_1)\subset Y_1;
\]

\item with $L:=F_x(\lambda_*,\beta_*,0)$,
\[
\ker L=\operatorname{span}\{x_1,x_2\},\qquad
\ker L\cap X_i=\operatorname{span}\{x_i\},
\]
and $L$ is Fredholm of index zero;

\item $LX_i\subset Y_i$ and
$Y_i\cap\mathcal R(L)$ has codimension one in $Y_i$.
Let $P_i:Y\to Y_i$ be the bounded projections associated with
$Y=Y_1\oplus Y_2$.  Choose nonzero continuous functionals
$\ell_i^0\in Y_i^*$ such that
\[
\ker\ell_i^0=Y_i\cap\mathcal R(L),
\]
and extend them to all of $Y$ by
\begin{equation}\label{eq6.extfun}
\ell_i:=\ell_i^0\circ P_i\in Y^*,\qquad i=1,2.
\end{equation}
Then $\ell_i$ annihilates $\mathcal R(L)$ and its normalization is immaterial;

\item
\[
A_i:=\ell_i(F_{\lambda x}[1,x_i])\neq0,
\qquad i=1,2;
\]

\item with
\[
C_i:=\ell_i(F_{\beta x}[1,x_i]),
\]
the parameter determinant satisfies
\[
A_1C_2-A_2C_1\neq0;
\]

\item writing
\[
B_{11}:=\ell_1(F_{xx}[x_1,x_1]),
\qquad
B_{12}:=\ell_2(F_{xx}[x_1,x_2]),
\]
one has
\begin{equation}\label{eq6.quaddet}
A_1B_{12}-\frac12A_2B_{11}\neq0.
\end{equation}
\end{enumerate}
All derivatives above are evaluated at $(\lambda_*,\beta_*,0)$.

Then, for $\beta=\beta_*+\xi$ with $0<|\xi|$ sufficiently small, a secondary
branch bifurcates from a nontrivial point of the primary branch in $X_1$.
More precisely, after a smooth reparametrization there exist smooth functions
$\chi,\kappa$ with $\kappa(0,0)\neq0$ and a smooth complement map $z$ such that
\[
\widetilde\lambda(\xi,s)
=\lambda_*+\xi\chi(\xi,s),
\]
and
\[
\widetilde x(\xi,s)
=
\xi\kappa(\xi,s)x_1+\xi s x_2
+z\!\left(
\widetilde\lambda(\xi,s),\beta_*+\xi,
\xi\kappa(\xi,s),\xi s
\right).
\]
At $s=0$ this branch lies on the nontrivial primary branch.
\end{theorem}

\begin{proof}
This is the Lyapunov--Schmidt form of the secondary-bifurcation argument used in
\cite{Shearer}; we include the reduction because only this specialized version is
needed here.

Choose closed complements $Z_i\subset X_i$ such that
\[
X_i=\operatorname{span}\{x_i\}\oplus Z_i,
\qquad i=1,2,
\]
and set $Z:=Z_1\oplus Z_2$.  From (B2)--(B3),
\[
L:Z_i\longrightarrow Y_i\cap\mathcal R(L)
\]
is an isomorphism for each $i$, and therefore
$L:Z\to\mathcal R(L)$ is an isomorphism.  We use the corresponding range projection
that respects the splitting $Y=Y_1\oplus Y_2$.  Projecting the equation onto
$\mathcal R(L)$ and applying the implicit function theorem yields a smooth complement
\[
z=z_1+z_2=z(\lambda,\beta,a,b)\in Z_1\oplus Z_2
\]
for solutions written as
\[
x=ax_1+bx_2+z.
\]
Because the range equation and the complement are adapted to the two invariant
splittings, when $b=0$ the choice $z_2=0$ solves the $Y_2$ range equation; uniqueness
from the implicit function theorem therefore gives
\begin{equation}\label{eq6.zinvariance}
z(\lambda,\beta,a,0)\in Z_1\subset X_1.
\end{equation}

The remaining equations are the two scalar Lyapunov--Schmidt equations
\[
G_i(\lambda,\beta,a,b)
:=
\ell_i\!\left(F(\lambda,\beta,ax_1+bx_2+z)\right)=0,
\qquad i=1,2,
\]
where the extended functionals \eqref{eq6.extfun} are defined on all of $Y$.
Because $F(\lambda,\beta,0)=0$, the first reduced equation restricted to $b=0$
has the factor $a$.  By \eqref{eq6.zinvariance} and the invariance
$F(\lambda,\beta,X_1)\subset Y_1$, one has
$G_2(\lambda,\beta,a,0)=0$; hence, by the smooth Hadamard lemma, the second reduced
equation has the factor $b$.  Thus, near the origin,
\[
G_1(\lambda,\beta,a,0)=a\,H_1(\lambda,\beta,a),
\qquad
G_2(\lambda,\beta,a,b)=b\,H_2(\lambda,\beta,a,b).
\]
At the double point,
\[
\partial_\lambda H_1=A_1,\qquad
\partial_aH_1=\frac12B_{11},
\]
and
\[
\partial_\lambda H_2=A_2,\qquad
\partial_aH_2=B_{12}.
\]
Hence \eqref{eq6.quaddet} is exactly the nonvanishing Jacobian determinant of
$(H_1,H_2)$ with respect to $(\lambda,a)$.  The implicit function theorem therefore
solves
\[
H_1=H_2=0
\]
for $(\lambda,a)$ as smooth functions of $(\beta,b)$.

Differentiating these two equations with respect to $\beta$ at the double point and
$b=0$ shows that the derivative of the primary amplitude is nonzero precisely when
\[
A_1C_2-A_2C_1\neq0.
\]
Thus (B5) yields
\[
a=(\beta-\beta_*)\kappa(\beta-\beta_*,b),
\qquad \kappa(0,0)\neq0.
\]
Putting $\beta-\beta_*=\xi$ and $b=\xi s$ gives the asserted parametrization.
For $s=0$ one has $b=0$, so the solutions lie in $X_1$ and form the primary branch;
for $s\neq0$ they leave $X_1$, giving the secondary branch.
\end{proof}

\section*{Acknowledgments}
The work was supported by the National Natural Science Foundation of China (No. 12301133). Additional support was provided by the Postdoctoral Science Foundation of China (Nos. 2023M741441 and 2024T170353) and the Jiangsu Education Department (No. 23KJB110007).

\section*{Data availability}
No datasets were generated or analysed in this study.

\section*{Competing interests}
The authors declare that they have no
competing interests.

\end{document}